\documentstyle[12pt,srcltx,color,amsmath,amssymb]{article}
\newtheorem{lemma}{Lemma}[section]
\def\U{{\mathbf{U}}}
\def\Z{{\mathbf{Z}}}
\def\W{{\mathbf{W}}}
\def\V{{\mathbf{V}}}
\def\H{{\mathbf{H}}}
\def\pen{{{\hbox{\rm\scriptsize pen}}}}
\def\reg{{{\hbox{\rm\scriptsize reg}}}}
\def\Id{{\hbox{\rm Id}}}
\def\Col{{\hbox{\rm Col}}}
\def\Row{{\hbox{\rm Row}}}
\def\dist{{\hbox{\rm dist}}}
\def\cS{{\cal S}}
\def\cU{{\cal U}}

\def\cA{{\cal A}}
\def\cZ{{\cal Z}}
\def\cM{{\cal M}}
\def\cE{{\cal E}}
\def\cW{{\cal W}}
\def\cV{{\cal V}}
\def\cH{{\cal H}}

\def\cE{{\cal E}}

\def\Ker{\hbox{\rm Ker}}
\def\Argmin{\mathop{\rm Argmin}}
\def\bbr{{\mathbf{R}}}
\def\cP{{\cal P}}
\def\bC{{\mathbf{C}}}
\def\dim{\hbox{\rm dim}}

\def\Opt{{\hbox{\rm Opt}}}
\def\Rank{\hbox{\rm Rank}}

\def\Tr{\mathop{\hbox{\rm Tr}}}
\def\lf{{\hbox{\tiny left}}}
\def\rg{{\hbox{\tiny right}}}
\newtheorem{proposition}{Proposition}[section]

\newtheorem{remark}{Remark}[section]

\newcommand{\be}{\begin{eqnarray}}
\newcommand{\ee}[1]{\label{#1}\end{eqnarray}}

\newcommand{\ese}{\end{eqnarray*}}
\newcommand{\bse}{\begin{eqnarray*}}
\newcommand{\epr}{\hfill\hbox{\hskip 4pt
                \vrule width 5pt height 6pt depth 1.5pt}\vspace{0.5cm}\par}
\begin{document}
\title{On unified view of nullspace-type conditions for recoveries associated with general sparsity structures}
\author{Anatoli Juditsky\thanks{LJK,
Universit\'e J. Fourier, B.P. 53, 38041 Grenoble
Cedex 9, France, {\tt Anatoli.Juditsky@imag.fr}} \and Fatma K{\i}l{\i}n\c{c} Karzan\thanks{Carnegie Mellon University, Pittsburgh, Pennsylvania 15213, USA, {\tt fkilinc@andrew.cmu.edu}} \and Arkadi Nemirovski\thanks{Georgia Institute
 of Technology, Atlanta, Georgia
30332, USA, {\tt nemirovs@isye.gatech.edu}\newline
Research of
the third author was supported by the ONR grant N000140811104 and NSF grant DMS 0914785.}}
\maketitle
\begin{abstract}
We discuss a general notion of ``sparsity structure'' and associated recoveries of a sparse signal from  its linear image of reduced dimension possibly corrupted with noise. Our approach allows for unified treatment of (a) the ``usual sparsity'' and ``usual $\ell_1$ recovery,'' (b) block-sparsity with possibly overlapping blocks and associated block-$\ell_1$ recovery, and (c) low-rank-oriented recovery by nuclear norm minimization. The proposed recovery routines are natural extensions  of the usual $\ell_1$ minimization used in Compressed Sensing. Specifically, we present nullspace-type sufficient conditions for the recovery  to be precise on sparse signals in the noiseless case. Then we derive error bounds for imperfect  (nearly sparse signal, presence of observation noise, etc.) recovery under these conditions. In all of these cases, we present efficiently verifiable sufficient conditions for the validity of the associated nullspace properties.
\end{abstract}
\section{Introduction}
We address the problem of recovering a {\sl representation} $w=Bx\in E$ of an unknown signal $x\in X$ via noisy observations
\[
y=Ax+\xi
\]
of $x$. Here $X,E$ are Euclidean spaces, $A:X\to \bbr^m$ and $B:X\to E$ are given linear {\sl sensing} and {\sl representation} maps, and $\xi$ is ``uncertain-but-bounded''  observation error satisfying $\phi(\xi)\leq\epsilon$ ($\phi(\cdot)$ is a given norm on $\bbr^m$, and $\epsilon$ is a given error bound). We consider, for instance, the recovering routine of the form
\[
y\mapsto \widehat{x}(y)\in\Argmin_{u\in X}\left\{\|Bu\|: \phi(Au-y)\leq\epsilon\right\}\mapsto \widehat{w}(y)=B\widehat{x}(y),
\]
and we want this recovery to behave well provided that $Bx$ is sparse in some prescribed sense. In this note, we
introduce a rather general notion of {\sl sparsity structure} on the representation space $E$ which, under some structural restriction on the norm $\|\cdot\|$, allows to point out ``nullspace type'' conditions for the recovery to be precise provided that $Bx$ is ``$s$-sparse'' with respect to our structure. It also allows for explicit error bounds for ``imperfect recovery'' (noisy observations, near $s$-sparsity instead of the exact one, etc.) The motivation behind this construction is to present a simple unified general framework, which allows, for instance, for a simple treatment of three important particular cases:
\begin{itemize}
\item recovering $s$-sparse in the usual sense (at most $s$ nonzero entries) signals via $\ell_1$ minimization (the corresponding nullspace property goes back to \cite{NP1,NP2})
\item recovering $s$-block-sparse signals via block-$\ell_1$ minimization, and
\item recovering matrices of low rank via nuclear norm minimization.
\end{itemize}
We present the respective sparsity structures and provide {\sl verifiable} sufficient conditions for the validity of associated nullspace properties (and thus -- for the validity of the corresponding recovery routines); the prototypes of our verifiable conditions can be found in \cite{JuN1,JuKKN2,JuN3,Judetal}.

\section{Problem description and recovery routines}
\subsection{Situation}\label{situation} Let $X$, $E$ be Euclidean spaces, $x\mapsto Ax:X\to\bbr^m$ be a linear {\sl sensing map}, and $x\mapsto Bx: X\to E$ be a linear  {\sl representation map}. We are interested in the problem as follows:
\begin{quote}
(!) {\sl Given a noisy observation
\begin{equation}\label{eqobs}
y=Ax+\xi
\end{equation}
of an unknown signal $x\in X$, we want to recover the representation $w=Bx$ of $x$.}
\end{quote}
\paragraph{Sparsity structure}
We will focus on the case when a priori information on $x$ is expressed in terms of properly defined {\sl sparsity} of the representation $Bx$ of $x$. To this end, we define a {\sl sparsity structure} on $E$, specifically, as follows:
\begin{quote}
{\sl Let $\|\cdot\|$ be a norm on $E$, $\|\cdot\|_*$ be the conjugate norm, and $\cP$ be a family of linear maps of $E$ into itself, such that
\begin{itemize}
\item[{\bf A.1.}] Every $P\in\cP$ is a projector: $P^2=P$.
\item[{\bf A.2.}] Every $P\in\cP$ is assigned a nonnegative weight $\nu(P)$ and a linear map $\overline{P}$ on $E$ such that $\overline{P}P=0$;
\item[{\bf A.3.}] Whenever $P\in\cP$ and  $f,g\in E$, one has
\begin{equation}\label{(!)}
\|P^*f+\overline{P}^*g\|_*\leq \max[\|f\|_ *,\|g\|_*],
\end{equation}
where for a linear map $Q$ acting from a Euclidean space $E$ into a Euclidean space $F$, $Q^*$ is the conjugate mapping acting from $F$ to $E$.
\end{itemize}}
A collection  of the just introduced entities satisfying the requirements A.1-3 will be referred to as a {\sl sparsity structure} on $E$.
\end{quote}
From now on, given a sparsity structure, we set for a nonnegative real $s$
\[
\cP_s=\{P\in\cP: \nu(P)\leq s\}.
\]
Given $s\geq0$, we call a vector $w\in E$ {\sl $s$-sparse}, if there exists $P\in\cP$ such that $\nu(P)\leq s$ and $Pw=w$. We call a vector $x\in X$ as $s$-sparse, if so is its representation $w=Bx$.

\par
We are about to present some instructive examples of the just outlined situation. Given a finite set $V$, we denote by $\bbr(V)$ the space of real vectors with entries indexed by elements of $V$ and equip this space with the standard inner product. For a subset $S\in V$, we set $\bbr(S)=\{v\in\bbr(V):v_i=0\,\forall i\not\in S\}$, and refer to $\bbr(S)$ as a {\sl coordinate subspace} of $\bbr(V)$. By $P_S$, $S\subset V$, we denote  the {\sl coordinate projector} -- the natural projector of $\bbr(V)$ onto $\bbr(S)$.
\subsection{Examples}
\paragraph{Example I.a: $\ell_1$ recovery}
In this example,
 \begin{itemize}
 \item $X=E=\bbr^n=\bbr(V:=\{1,...,n\})$ with the standard inner product;
 \item The representation map is the identity: $Bx\equiv x$;
 \item $\cP$ is comprised of projectors onto all coordinate subspaces of $\bbr^n$, $\overline{P}=I_n-P$;
 \item $\nu(P)=\Rank(P)$;
 \item $\|\cdot\|=\|\cdot\|_1$.
 \end{itemize}
 Properties {\bf A.1-3} clearly take place, $s$-sparsity of $x\in X$ as defined above is the usual sparsity (at most $s$ nonzero entries), and (\ref{recovery})  is the usual $\ell_1$ recovery.
\paragraph{Example I.b: Group $\ell_1$ recovery}
Now we want to model the ``block-sparse'' situation as follows. $X=\bbr^n=\bbr(V=\{1,...,n\})$, and the index set $V$ is represented as the union of $K$ nonempty (and possibly overlapping) subsets $V_1,...,V_K$, so that to every $x\in X$ one can associate {\sl blocks} $x^\ell$, $\ell=1,...,K$, which are natural projections of $x$ onto  $\bbr(V_\ell)$. Assuming the subsets $V_\ell$ to be assigned with positive  weights $\chi_\ell$, we define the sparsity of a signal $x\in X$ as the ``weighted number of nonzero blocks $x^\ell$,'' that is, the quantity $s(x)=\sum_{\ell:x^\ell\neq0}\chi_\ell$.  This can be modeled by the following representation and sparsity structure:
\begin{itemize}
\item $E=\bbr(V_1)\times...\times \bbr(V_K)$, so that $w\in E$ is a block vector $[w^1;...;w^K]$ with $w^\ell\in\bbr(V_\ell)$, and $Bx=[x^1;...;x^K]\in E$;
\item $\cP$ is comprised of orthoprojectors $P=P_I$ onto the subspaces $E_I=\{[w^1;...;w^K]\in E: w^\ell =0\,\forall\ell\not\in I\}$, associated with subsets $I$ of the index set $\{1,...,K\}$, and $\overline{P}=\Id-P$, where $\Id$ stands for the identity mapping on $X$;
\item $\nu(P_I)=\sum_{\ell\in I} \chi_\ell$;
 \item The norm $\|\cdot\|$ is defined as follows. For every $\ell\leq K$, we equip $\bbr(V_\ell)$ with a norm $\|\cdot\|_{(\ell)}$, and set
\[
\|[z^1;...;z^K]\|=\sum_{\ell=1}^K\|z^\ell\|_{(\ell)}
\]
\end{itemize}
Verification of A.1-A.3 is immediate. Note that when $V_1,...,V_K$ do not overlap and $\chi_\ell=1$ for all $\ell$, we find ourselves in the standard block-sparse situation: $E$ can be naturally identified with $X$, making $B$ the identity, vectors from $X=\bbr^n$ are split into $K$ non-overlapping blocks, the norm is the block-$\ell_1$ norm, and $s$-sparsity of $x$ means that $x$ has at most $s$ nonzero blocks.
\paragraph{Example II: Nuclear norm recovery}
In this example,
\begin{itemize}
\item $X=E=\bbr^{p\times q}$ with $p\geq q$ and the Frobenius inner product, $B$ is the identity;
\item $\cP$ is comprised of the mappings $P(x)=P_\lf x P_\rg$, where $P_\lf\in\bbr^{p\times p}$ and $P_\rg\in\bbr^{q\times q}$ are orthoprojectors, and $\overline{P}(x)=(I_p-P_\lf)x(I_q-P_\rg)$;\footnote{To avoid notational ambiguity, in the situation of Example II we denote the image of $z\in E$ under a mapping $P\in\cP$ by $P(z)$ rather than by $Pz$, to avoid collision with the notation for matrix products like $(I_p-P_\lf)x(I_q-P_\rg)$.}
\item $\|\cdot\|$ is the nuclear norm: $\|x\|=\sum_{j=1}^q \sigma_j(x)$, where $\sigma_1(x)\geq\sigma_2(x)...\geq\sigma_q(x)$ are the singular values of $x\in\bbr^{p\times q}$.
\end{itemize}
The assumptions {\bf A.1-3} clearly take place.
\begin{quote}
{\small
To verify {\bf A.3}, observe that the norm conjugate to $\|\cdot\|$ is the spectral norm $\|x\|_{2,2}:=\sigma_1(x)$. We now have
\bse
&&\|f\|_{2,2}\leq 1,~~\|g\|_{2,2}\leq 1\\
&\Rightarrow& \|P^*(f)\|_{2,2}=\|P_\lf f P_\rg\|_{2,2}\leq 1,~~ \|\overline{P}^*(g)\|_{2,2}=\|(I_p-P_\lf)g(I_q-P_\rg)\|_{2,2}\leq1\\
&\Rightarrow& \|P^*(f)+\overline{P}^*(g)\|_{2,2}\leq1,
\ese
where the last $\Rightarrow$ follows from the fact that the orthogonal complements to the kernels of $P^*(f)$ and $\overline{P}^*(g)$, same as the images of these matrices, are orthogonal to each other.}
\end{quote}
In this case, by assigning the projectors $P\in\cP$ with weights according to
\[
\nu(P)=\max[\Rank(P_\lf),\Rank(P_\rg)],
\]
we arrive at the situation where $s$-sparse signals  are exactly the $p\times q$ matrices of rank $\leq s$.
\section{Main results}
 \subsection{Recovery routines} Let $\phi(\cdot)$ be a norm on the image space of $A$, and let $\epsilon$ be an a priori upper bound on the $\phi$-norm $\phi(\xi)$ of the observation error, see (\ref{eqobs}). In order to recover the representation $Bx$ of signal $x$ underlying the observation (\ref{eqobs}), we use {\sl regular} recovery -- the ``standard'' recovery by  $\|\cdot\|$-minimization as follows:
\begin{equation}\label{recovery}
y\mapsto\widehat{x}_\reg(y)\in\Argmin_{u\in X}\left\{\|Bu\|: \phi(Au-y)\leq\epsilon\right\}\mapsto \widehat{w}_\reg(y)=B\widehat{x}_\reg(y)
\end{equation}
and we treat $\widehat{w}_\reg(y)$ as the resulting estimate of $Bx$.
\begin{quotation}
{\em We say that the sensing map $A$ is {\em $s$-good} if the above recovery {\sl in the noiseless case ($\epsilon=0$)} reproduces exactly the representation $Bx$ of every $s$-sparse signal $x$.}
\end{quotation}
We also consider an alternative to (\ref{recovery}), specifically, the {\sl penalized $\|\cdot\|$ recovery} introduced in \cite{JuN3,Judetal}. This recovery routine is given by
\begin{equation}\label{penalized}
y\mapsto \widehat{x}_\pen(y)\in\Argmin_{u\in X}\left\{\|Bu\|+\lambda\phi(Au-y)\right\}\mapsto \widehat{w}_\pen(y)=B\widehat{x}_\pen(y)
\end{equation}
where $\lambda>0$ is a penalty parameter.
\subsection{$s$-goodness and nullspace property}\label{immediate}
We start with the following immediate observation:
\begin{lemma}\label{lem1} In the situation of section \ref{situation}, let $w\in E$, $P\in \cP$ be such that $Pw=w$. Then for every $z\in X$ one has
\begin{equation}\label{eq1}
\|w+Bz\|\geq \|w\|+\|\overline{P}Bz\|-\|PBz\|.
\end{equation}
In particular, if the following ``nullspace property'' holds true:
\begin{equation}\label{eq2}
\forall P\in \cP_s,\;z\in\Ker(A),Bz\neq0:~ \|\overline{P}Bz\|>\|PBz\|,
\end{equation}
then $A$ is $s$-good.
\end{lemma}
{\bf Proof.} Let $w\in E$ and $P\in \cP$ be such that $Pw=w$. We have
\[
\forall (f,g:\max\{\|f\|_*,\|g_*\}\leq1): ~~\|P^*f+\overline{P}^*g\|_*\leq 1~~~\hbox{\ [by A.3]}
\]
which implies
\bse
& \Rightarrow\|w+Bz\| &\geq \|w+Bz\| \|P^*f+\overline{P}^*g\|_* \\
&& \geq \langle P^*f+\overline{P}^*g,w+Bz\rangle=\langle f,Pw+PBz\rangle +\langle g,\overline{P}w+\overline{P}Bz\rangle\\
&&\geq \langle f,w\rangle +\langle f,PBz\rangle + \langle g,\overline{P}Bz\rangle ~~~\hbox{\ [since $w=Pw$ and therefore $\overline{P}w=\overline{P}Pw=0$]}\\
&&\geq \langle f,w\rangle -\|PBz\|+\langle g,\overline{P}Bz\rangle.
 \ese
 When choosing $f,g$ to be such that $\|f\|_*=\|g\|_*=1$ and $\langle f,w\rangle=\|w\|$, $\langle g,\overline{P}Bz\rangle =\|\overline{P}Bz\|$, we get
 \[
 \|w+Bz\|\geq \|w\|-\|PBz\|+\|\overline{P}Bz\|,
 \]
 as claimed in (\ref{eq1}). \par
 Now let (\ref{eq2}) take place, and let us prove that $A$ is $s$-good. All we need to verify is that if $x$ is $s$-sparse and $\widehat{x}=\widehat{x}(Ax)$, see (\ref{recovery}), then $Bx=B\widehat{x}$. Setting $z=\widehat{x}-x$, so that $z\in\Ker (A)$, choosing $P\in\cP$ such that $\nu(P)\leq s$ and $PBx=Bx$ and applying (\ref{eq1}) with $w=Bx$, we get
 \[
 \|B\widehat{x}\|=\|Bx+Bz\|\geq \|Bx\|+\|\overline{P}Bz\|-\|PBz\|,
 \]
 while by the origin of $\widehat{x}$ we have $\|B\widehat{x}\|\leq \|Bx\|$. It follows that $\|\overline{P}Bz\|\leq \|PBz\|$, which, by the nullspace property (\ref{eq2}), is possible only when $Bz=0$ (recall that $\nu(P)\leq s$ and $z\in\Ker (A)$). \epr
\begin{remark}\label{rem1}{\rm Independently of any assumptions on $\|\cdot\|$, an evident {\sl necessary} condition for $s$-goodness of $A$ is:
 \begin{quote}
 {\sl Whenever $P\in\cP$, $\nu(P)\leq s$, $x\in X$ are such that $PBx=Bx$, and $z\in\Ker(A)$, $Bz\neq0$,  there exists $f$, $\|f\|_*=1$, such that $\langle f,PBx\rangle=\|PBx\|$ and $\langle f,Bz\rangle \geq0$},
 \end{quote}
 When modifying this condition by replacing $\langle f,Bz\rangle\geq0$ in the conclusion with $\langle f,Bz\rangle>0$, this necessary condition for $s$-goodness becomes sufficient. This, by the way, immediately implies the necessity and sufficiency of the standard nullspace property \cite{NP1,NP2} for the validity of $\ell_1$ minimization and the translation of it to the matrix case given by \cite{Tuncel,OymakFazel2011}.
Moreover, recently, \cite{CandesRecht2012} has established the sufficiency of this condition in the case of decomposable norms, i.e., sparse, non-overlapping block-sparse and low-rank recovery following a unified view based on subdifferentials.}
 \end{remark}

 \subsection{Error bounds for imperfect $\|\cdot\|$ recovery}\label{simperfect}

 \paragraph{Conditions $\bC_s(\gamma,\beta;\phi)$ and $\bC^+_s(\gamma,\beta;\phi)$.} In the sequel, we shall use the following two conditions (where $s,\gamma\geq0$, $\beta\in[0,\infty]$, and $\phi(\cdot)$ is a norm on the image space of $A$):
  \begin{quote}
 $\bC_s(\gamma,\beta;\phi)$:
 \[
 \forall (z\in X,P\in \cP_s): ~~\|PBz\|+\|Bz\|-\|\overline{P}Bz\|\leq \beta\phi(Az)+\gamma\|Bz\|
 \]
 \end{quote}
 (from now on, $(+\infty)\cdot 0=0$ and $(+\infty)\cdot a=+\infty$ when $a>0$).
 \begin{quote}
 $\bC^+_s(\gamma,\beta;\phi)$: there exists a (semi)norm $\|\cdot\|_{\cP_s}$ on $E$ such that
 \begin{equation}\label{bCplus}
 \begin{array}{ll}
 (a)& \forall (z\in X,P\in \cP_s):~~ \|PBz\|+\|Bz\|-\|\overline{P}Bz\|\leq \|Bz\|_{\cP_s}\\
 (b)&\forall z\in X:~~ \|Bz\|_{\cP_s}\leq\beta\phi(Az)+\gamma\|Bz\|.\\
 \end{array}
 \end{equation}
 \end{quote}
 Let us make the following immediate observations:
 \begin{remark}\label{rem3} {\rm (i)} The validity of condition $\bC_s(\gamma,\beta;\phi)$ with some $\gamma<1$, and some $\beta$, $\phi$ implies the validity of the nullspace property {\rm (\ref{eq2})};
 \par
 {\rm (ii)} $\bC^+_s(\gamma,\beta;\phi)$ implies $ \bC_s(\gamma,\beta;\phi)$;
 \par
 {\rm (iii)} The (semi)norm $\|w\|_{\cP_s}^+:=\max\limits_{P\in\cP_s}\left[\|Pw\|+\|(\Id-\overline{P})w\|\right]$ satisfies {\rm (\ref{bCplus}.$a$)} {\rm (since by Triangle inequality we have $\|PBz\|+\|Bz\|-\|\overline{P}Bz\|\leq \|PBz\|+\|(\Id-\overline{P})Bz\|$)};
 \par
 {\rm (iv)} Whenever $s'\leq s,\;\gamma'\geq \gamma$, and $\beta'\geq\beta$ one has  $\bC^+_{s}(\gamma,\beta;\phi)\Rightarrow \bC^+_{s'}(\gamma',\beta';\phi)$, and similarly for the condition $\bC_s$;\par
 {\rm (iv)} Condition $\bC^+_s(\gamma,+\infty;\phi)$ reads
 \begin{quote} There exists a (semi)norm $\|\cdot\|_{\cP_s}$ on $E$ such that
\bse
& (a)& \forall (z\in X,P\in \cP_s):~~ \|PBz\|+\|Bz\|-\|\overline{P}Bz\|\leq \|Bz\|_{\cP_s}\\
& (b)&\forall z\in\Ker(A):~~ \|Bz\|_{\cP_s}\leq\gamma\|Bz\|.
 \ese
 \end{quote}
 if this condition is satisfied, then for every $\gamma'>\gamma$ there exists $\beta<\infty$ such that the condition $\bC^+_s(\gamma',\beta;\phi)$ is satisfied.
 \end{remark}
 \paragraph{Error bounds for imperfect regular $\|\cdot\|$ recovery} are stated in the following
 \begin{proposition}\label{proprec} In the situation of section \ref{situation}, let a sparsity level $s\geq0$ be given, and let the condition $\bC_s(\gamma,\beta;\phi)$ take place for some $\gamma<1$ and $\beta<\infty$. Given a signal $x\in X$ along with its observation $y=Ax+\xi$, where $\phi(\xi)\leq\epsilon$, let $x$ be ``nearly $s$-sparse'' and $\widehat{x}$ be ``nearly $\widehat{x}_\reg(y)$,'' specifically, for given nonnegative tolerances $\delta_x$, $\delta_\phi$, and $\delta$ one has
 \begin{itemize}
 \item  there exists $P\in\cP_s$ such that $\|(I-P)Bx\|\leq \delta_x$ (``near $s$-sparsity of $x$'');
 \item one has
 \begin{equation}\label{meaning}
 \begin{array}{ll}
 (a)&\phi(A\widehat{x}-[Ax+\xi])\leq\epsilon+\delta_\phi\\
 (b)&\|B\widehat{x}\|\leq \Opt+\delta,\,\,~\Opt:=\min\limits_u\left\{\|Bu\|:\phi(Au-y)\leq\epsilon\right\}.\\
 \end{array}
 \end{equation}
 (``$\widehat{x}$ is a nearly feasible nearly optimal solution to the optimization problem specifying $\widehat{x}(y)$'').
 \end{itemize}
 Then
 \begin{equation}\label{errorbound}
 \|B\widehat{x}-Bx\|\leq {\beta[2\epsilon+\delta_\phi]+\delta+2\delta_x\over 1-\gamma}.
 \end{equation}
 \end{proposition}
 Proofs of the results of this section are put in the appendix.

\paragraph{Error bounds for penalized $\|\cdot\|$ recovery}%
\begin{proposition}\label{penrec} In the situation of section \ref{situation}, let $s\geq0$, $\gamma<1$, $\phi(\cdot)$ and $\beta<\infty$ be such that the condition $\bC_s(\gamma,\beta;\phi)$ is satisfied. Let also the penalty parameter $\lambda$ in {\rm (\ref{penalized})} be $\geq\beta$. Finally, let the signal $x$ underlying the observations  be ``nearly $s$-sparse,'' meaning that there exists $P\in\cP_s$ such that $\|(I-P)Bx\|\leq\delta_x$, and let $\widehat{x}$ be a near-optimal solution to the optimization problem specifying $\widehat{x}_\pen$, namely,
$$
\lambda \phi(A\widehat{x}-y) +\|B\widehat{x}\|\leq \min_z\left\{\lambda\phi(Az-y) +\|Bz\|\right\}+\delta.
$$
Then
\begin{equation}\label{imperfect}
\|B\widehat{x}-Bx\|\leq {2\delta_x+\delta+2\lambda\phi(\xi)\over1-\gamma}
\end{equation}
where $\xi=y-Ax$ is the observation noise.
\end{proposition}

\paragraph{Comment:} As compared to the plain $\|\cdot\|$-recovery (\ref{recovery}), the penalized $\|\cdot\|$ recovery requires a priori knowledge of a $\beta<\infty$ such that the condition $\bC_s(\gamma,\beta;\phi)$  takes place; indeed, in order for the error bound (\ref{imperfect}) to be applicable, we should ensure $\lambda\geq\beta$. As a compensation, the penalized $\|\cdot\|$ recovery does not require any a priori information on the level of observation error, and as such
is well suited for the case when the latter is random (or a sum of a random and a bounded components).

\section{Application examples}
In the rest of this note, we are interested in the particular forms taken by the nullspace sufficient condition for $s$-goodness of $A$ (Lemma \ref{lem1}) and the error bound for imperfect $\|\cdot\|$-recovery (Proposition \ref{proprec}) in the examples described in section \ref{situation}.
\subsection{Example I.a: $\ell_1$ recovery}
In the situation of Example I.a,  the nullspace property (\ref{eq2}) reads
 \begin{equation}\label{Nullspaceell1}
 \gamma_s(A):=\max_x\left\{\|x\|_{s,1}:x\in\Ker (A),\|x\|_1\leq 1\right\}<1/2,
 \end{equation}
 where $\|x\|_{s,1}$ is the sum of the $s$ largest magnitudes of entries in $x$. This is a well-known necessary and sufficient condition for the validity of the standard sparse $\ell_1$ recovery \cite{NP1,NP2}. It is immediately seen that condition $\bC_s^+(\gamma,\beta;\phi)$ is satisfied if and only if it is satisfied when setting $\|z\|_{\cP_s}=2\|z\|_{s,1}$, and in this case the condition is equivalent to $\bC_s(\gamma,\beta;\phi)$. The latter condition reads
 \begin{equation}\label{conditionCell1}
 \forall (z\in\bbr^n): ~~\|z\|_{s,1}\leq {\beta\over 2}\phi(Az)+{\gamma\over 2}\|z\|_1.
 \end{equation}
 Validity of this relation is equivalent to the fact that the quantity $\widehat{\gamma}_s(A,\cdot)$ introduced in \cite{JuN1} satisfies $\widehat{\gamma}_s(A,\beta/2)\leq \gamma/2$ (see \cite[Theorem 2.2]{JuN1}.  What is denoted by $\|\cdot\|$ in the latter reference, is now $\phi(\cdot)$), and with this in mind, error bounds (\ref{errorbound}) and (\ref{imperfect}) recover the bounds in \cite[Theorem 3.1]{JuN1}. Beside this, one can find in \cite{JuN1} verifiable sufficient conditions for the validity of $\bC_s(\gamma,\beta;\phi)$, their relations to Restricted Isometry Property, etc.
\subsection{Example I.b: Group $\ell_1$ recovery}
In the situation of Example I.b, given a positive real $s$, let us define  the norm $\pi_s(\cdot)$ on $\bbr^K$ as
\bse
\pi_s(u)=2\max\limits_\eta\left\{\sum_{\ell=1}^K\eta_\ell |u_\ell|:\eta_\ell\in\{0,1\}\,\forall \ell,\,~\sum_{\ell=1}^K \chi_\ell\eta_\ell\leq s\right\},
\ese
and let $\|\cdot\|_{1,s}$ be the norm on $E$ given by
\[
\|w\|_{1,s}=\pi_s([\|w^1\|_{(1)};\|w^2\|_{(2)};...;\|w^K\|_{(K)}]).
\]
Observe that for every $z\in X$ and every $P=P_I\in\cP_s$ we have
\bse
\|PBz\|+\|Bz\|-\|\overline{P}Bz\| &=&\sum_{\ell\in I}\|(Bz)^\ell\|_{(\ell)}+\sum_{\ell=1}^{K} \|(Bz)^\ell\|_{(\ell)}-\sum_{\ell\not\in I}\|(Bz)^\ell\|_{(\ell)}\\
&=&2\sum_{\ell\in I}\|(Bz)^\ell\|_{(\ell)}\leq \|Bz\|_{1,s},
\ese
where the concluding $\leq$ is given by $P_I\in\cP_s$, so that $\sum_{\ell\in I} \chi_\ell\leq s$. Thus, with $\|\cdot\|_{\cP_s}$ set to $\|\cdot\|_{1,s}$ the condition (\ref{bCplus}) is satisfied. We have arrived at the following
\begin{proposition}\label{propnover} In the situation of Example I.b, for every positive reals $s$, $\beta,\gamma$, the condition
\begin{equation}\label{eq256}
\forall z\in X: ~~\|Bz\|_{1,s}\leq \beta\phi(Az)+\gamma\|Bz\|
\end{equation}
is sufficient for the validity of
$\bC^+_s(\gamma,\beta;\phi)$, and thus --- for the validity of $\bC_s(\gamma,\beta;\phi)$. As a result, condition {\rm (\ref{eq256})} with $\gamma<1$  is sufficient for $s$-goodness of $A$ and for the applicability of the error bounds {\rm (\ref{errorbound})} and (\ref{imperfect}).
\end{proposition}
\paragraph{A verifiable sufficient condition for (\ref{eq256})} Condition (\ref{eq256}) is difficult to verify. We are about to present a verifiable sufficient condition for the validity of (\ref{eq256}) inspired by \cite{Judetal}.  For a linear map $Q^{k\ell}:\bbr(V_\ell)\to \bbr(V_k)$, let $\|Q^{k\ell}\|_{(\ell k)}$ be the norm of the map induced by the norms $\|\cdot\|_{(\ell)}$ and $\|\cdot\|_{(k)}$ on the argument and the image spaces:
$$
\|Q^{k\ell}\|_{(\ell k)}=\max_{u}\{\|Q^{k\ell}u\|_{(k)}:\|u\|_{(\ell)}\leq1\}.
$$
Let also $N$ be the dimension of $E$, and let $n_1,...,n_K$ be the cardinalities of $V_1,...,V_K$, so that $E$ can be identified with $\bbr^N=\bbr^{n_1}\times...\times\bbr^{n_K}$. For an $N\times N$ matrix $W$, let $W^{k\ell}$ be the $n_k\times n_\ell$ blocks of $W$ associated with the direct product representation $\bbr^{n_1}\times...\times\bbr^{n_K}$ of $\bbr^N$, and let $\Omega[W]$ be the $K\times K$ matrix with the entries $\Omega_{k\ell}=\|W^{k\ell}\|_{(\ell k)}$, $1\leq k,\ell \leq K$.
\begin{proposition}\label{propsuffcond} In the situation of Example I.b, given $\gamma>0$, let $m\times N$ matrix $H$ and $N\times N$ matrix $W$ satisfy the relations
\begin{equation}\label{onehas}
\begin{array}{ll}
(a)&B=WB+H^TA\\
(b)&\max_{\ell\leq K}\pi_s\left(\Col_\ell(\Omega[W])\right)\leq\gamma,\\
\end{array}
\end{equation}
where $\Col_j(Q)$ is the $j$-th column of matrix $Q$.
Then the relation {\rm (\ref{eq256})} holds true with
\bse
\beta=\beta[H]:=\max_{v\in\bbr^m}\left\{\|H^Tv\|_{1,s}:~\phi(v)\leq1\right\}.
\ese
\end{proposition}
{\bf Proof.} Let $z\in X$. Under the premise of the proposition, we have
\bse
(Bz)^k&=&(WBz)^k+(H^TAz)^k=\sum_{\ell=1}^KW^{k\ell}(Bz)^\ell+(H^TAz)^k\\
\Rightarrow\|(Bz)^k\|_{(k)}&\leq& \sum_{\ell=1}^K \|W^{k\ell}\|_{(\ell k)}\|(Bz)^\ell\|_{(\ell)}+\|(H^TAz)^k\|_{(k)}\\
\Rightarrow [\|(Bz)^1\|_{(1)};...;\|(Bz)^K\|_{(K)}]&\leq& \sum_{\ell=1}^K\|(Bz)^\ell\|_{(\ell)}\Col_\ell(\Omega[W])+[\|(H^TAz)^1\|_{(1)};...;\|(H^TAz)^K\|_{(K)}]\\
\Rightarrow \|Bz\|_{1,s}&\leq&\|Bz\|\max_{\ell\leq K}\pi_s\left(\Col_\ell(\Omega[W])\right)+\|H^TAz\|_{1,s}
\ese
and the desired conclusion follows. \epr
\paragraph{Discussion} The sufficient condition for the validity of (\ref{eq256}) stated in Proposition \ref{propsuffcond} reduces to solving a system of convex constrains in matrix variables $H,W$ and scalars $\gamma$, $\beta$, namely,
\be
B=WB+H^TA,\; ~\pi_s\left(\Col_\ell(\Omega[W])\right)\leq{\gamma}\;\forall \ell,\;~\Psi_s(H):=\max\limits_{v:\phi(v)\leq1} \|H^Tv\|_{1,s}\leq \beta.
\ee{calS}
This system, although convex, still can be difficult to process, since the convex functions $\|W^{k\ell}\|_{(\ell k)}$, $\pi_s(\cdot)$ and $\Psi_s(H)$ can be difficult to compute. In such a case, one can replace these functions with their efficiently computable upper bounds (for details, ``solvable cases,'' etc., in the case of $\chi_\ell\equiv 1$ see
\cite{Judetal}). For instance,  (\ref{calS}) is computationally tractable when
\begin{itemize}
\item all norms $\|\cdot\|_{(\ell)}$ are either $\ell_1$, or $\ell_\infty$, or $\ell_2$ norms (this makes the matrix $\Omega[W]$ efficiently computable);
\item in appropriate scale, all weights $\chi_\ell$ are integers from a once for ever fixed (or polynomially growing with problem's sizes) range, which makes the norm $\pi_s$ efficiently computable. Note that one can always replace $\pi_s$ with a reasonably tight upper bound on $\pi_s$, specifically, the norm
    $$
    \widehat{\pi}_s(u)=2\max\limits_\eta\left\{\sum_{\ell=1}^K\eta_\ell |u_\ell|:0\leq\eta_\ell\leq\min[1,\hbox{Floor}(s/\chi_\ell)]\,\forall \ell,\,\sum_{\ell=1}^K \chi_\ell\eta_\ell\leq s\right\};
    $$
\item $\phi(\cdot)$ is the $\ell_1$ norm.
\end{itemize}
The last assumption is indeed restrictive. It, however, is responsible solely for the tractability of the constraint    $\Psi_s(H):=\max\limits_{v:\phi(v)\leq1} \|H^Tv\|_{1,s}\leq{\beta}$, and is crucial only if one's objective is to compute an upper bound on $\beta$.  On the other hand, it is of primary importance to ensure $\gamma<1$, since otherwise Proposition \ref{proprec} provides no error bound at all.
\par
Finally, we refer the reader to \cite{Judetal} for details on relationship of the derived conditions with block-RIP and other sufficient conditions used in the Compressive Sensing literature.
\subsection{Example II: Nuclear norm recovery}
For a $p\times q$ matrix $z$ with $p\geq q$, let
\[
\Sigma_k(z)=\sum_{i=1}^k\sigma_i(z),\,\,1\leq k\leq q.
\]
Observe that in the situation of Example II, where $B$ is the identity, and the sparsity parameter $s$ can be w.l.o.g.\@ restricted to be a nonnegative integer, we have (everywhere in this section, $\|\cdot\|$ is the nuclear norm).
\[
\forall (z\in\bbr^{p\times q},P\in\cP_s):~~\|P(z)\|\leq\Sigma_s(z),\,~\|\overline{P}(z)\|\geq\|z\|-\Sigma_{2s}(z).
\]
\begin{quote}
{\small Indeed, let $P\in\cP_s$, so that $\Rank(P_\lf)\leq s$ and $\Rank(P_\rg)\leq s$. Then $\|P(z)\|=\|P_\lf z P_\rg\|\leq \Sigma_s(z)$ by the Singular Value Interlacement Theorem. Since the matrix $\overline{P}(z)$ differs from $z$ by matrix of rank at most $2s$, by the same Singular Values Interlacement Theorem we have $\sigma_i(\overline{P}(z))\geq \sigma_{i+2s}(z)$, whence $\|\overline{P}(z)\|\geq \|z\|-\Sigma_{2s}(z)$.}
\end{quote}
We have arrived at the following
\begin{proposition}\label{proprank} In the situation of Example II, the norm
\[
\|z\|_{\cP_s}:=\Sigma_s(z)+\Sigma_{2s}(z)
\]
on $X=E=\bbr^{p\times q}$
satisfies the condition {\rm (\ref{bCplus}.$a$)}, so that the condition
\begin{equation}\label{elln}
\forall z\in \bbr^{p\times q}: ~\Sigma_s(z)+\Sigma_{2s}(z)\leq\beta\phi(Az)+\gamma\|z\|
\end{equation}
is sufficient for the validity of $\bC^+_s(\gamma,\beta;\phi)$, and thus --- for the validity of $\bC_s(\gamma,\beta;\phi)$. As a result, condition {\rm (\ref{elln})} with $\gamma<1$  is sufficient for $s$-goodness of $A$ and for applicability of the error bounds {\rm (\ref{errorbound})} and {\rm (\ref{imperfect})}.
\end{proposition}
Clearly there is a gap between the above sufficient condition and the necessary nullspace condition for exact low-rank matrix recovery, which is
 \[
2\Sigma_{s}(z) < \|z\|~~~\hbox{ for all } z\in\Ker(A), ~z\neq0.
\]
On the other hand our sufficient condition is stronger than, the only known to us, sufficient condition given in \cite{OymakFazel2011}, which requires
\[
2\Sigma_{2s}(z) < \|z\|~~~\hbox{ for all } z\in\Ker(A), ~z\neq0.
\]
\paragraph{A verifiable sufficient condition for (\ref{elln})} Following the same exposition scheme as in Examples I.a-b, it is now time to point out a verifiable sufficient condition for the validity of (\ref{elln}).
Let $H$ be a linear map from $X=\bbr^{p\times q}$ into $\bbr^m$, so that $W=\Id-H^*A$, $\Id$ being the identity mapping on $X$, is a linear map from $X$ into $X$. Assume that $W$ satisfies the requirement
\begin{equation}\label{eq500}
\forall z\in X: ~~\Sigma_s(Wz)+\Sigma_{2s}(Wz)\leq \gamma\|z\|;
\end{equation}
we claim that then (\ref{elln}) holds true with
\begin{equation}\label{eq501}
\beta=\beta[H]=\max\limits_{v} \left\{\Sigma_s(H^*v)+\Sigma_{2s}(H^*v):~\phi(v)\leq1\right\}.
\end{equation}
\begin{quote}
{\small Indeed, let $z\in X$. Assuming (\ref{eq500}), and setting $\pi(w)=\Sigma_s(w)+\Sigma_{2s}(w)$, $w\in\bbr^{p\times q}$, we have
$$
\begin{array}{l}
\pi(z)=\pi(Wz+H^*Az)\leq \pi(Wz)+\pi(H^*Az)\leq \gamma \|z\|+\beta[H]\phi(Az)\\
\end{array}
$$
as required in (\ref{elln}).}
\end{quote}
The question is, how to  efficiently certify the validity of (\ref{eq500}). The proposed answer is as follows.
Note that
\[
\Sigma_k(w)=\max_h\{\Tr(wh^T): \|h\|\leq k,\|h\|_{2,2}\leq1\},
\]
therefore (\ref{eq500}) is exactly
\bse
\gamma\geq\Opt[W]&:=&\max_{z\in\bbr^{p\times q}}\left\{\Sigma_s(Wz)+\Sigma_{2s}(Wz):~\|z\|\leq1\right\}\\
&=&\max\limits_{u,v,z\in\bbr^{p\times q}} \left\{\Tr([Wz][u+v]^T):~ \|z\|\leq 1,\;\begin{array}{l}\|u\|\leq s,\|u\|_{2,2}\leq1,\\
\|v\|\leq2s,\,\|v\|_{2,2}\leq1
\end{array}
\right\}.
\ese
Now,
$$
\Tr([Wz]h^T)=\sum_{i,j=1}^{pq} (\Theta[W])_{ij}(h^T\otimes z)_{ij},
$$
where $\Theta[W]$ is a properly defined {\sl linear in $W$} $pq\times pq$ matrix, and $\otimes$ is the Kronecker product; in other words $h^T\otimes z$ is the
$pq\times pq$ matrix with $p\times q$ blocks $[h^T\otimes z]^{\mu\nu}=h_{\nu\mu}z$, $1\leq \mu\leq q$, $1\leq \nu\leq p$.
We conclude that
\[
\Opt[W]\leq\max_{\cU,\cV}\left\{\sum_{i,j=1}^{pq}(\Theta[W])_{ij}[\cU_{ij}+\cV_{ij}]: \;\cU\in \Z_{s},\,\cV\in \Z_{2s}\right\},
\]
where
$$\Z_k=\left\{h^T\otimes w\in\bbr^{pq\times pq}\left|\begin{array}{l}h\in\H_k=\{h\in \bbr^{p\times q}:\|h\|\leq k,~\|h\|_{2,2}\leq1\}\\
w\in\W=\{w\in\bbr^{p\times q}:\|w\|\leq 1\}\\
\end{array}\right.\right\}.$$
When replacing $\Z_k$ with a computationally tractable set $\Z_k^*\supseteq\Z_k$ and setting
\begin{equation}\label{overline}
\overline{\Opt}[W]=\max_{\cU,\cV}\left\{\sum_{i,j=1}^{pq}(\Theta[W])_{ij}[\cU_{ij}+\cV_{ij}]: \;\cU\in \Z^*_{s},\,\cV\in \Z^*_{2s}\right\},
\end{equation}
we obtain an efficiently computable convex in $W$ upper bound on $\Opt[W]$, so that the efficiently computable convex constraint
\begin{equation}\label{eccc}
\overline{\Opt}[W]\leq\gamma
\end{equation}
implies the validity of (\ref{eq500}). We are about to point out two computationally tractable convex relaxations $\Z^*_k\supset
\Z_k$ of the sets $\Z_k$.

\par
Our first observation is as follows: since
the singular values of $h^T\otimes z$ are pairwise products of singular values of $h$ and  $z$, the set
\[
\U_k=\{\cU\in \bbr^{pq\times pq}:\; \|\cU\|\le k,\,\|\cU\|_{2,2}\le 1\}
\]
contains $\Z_k$, so that we can set $\Z^*_k=\U_k$. As a result, the efficiently computable convex constraint
\be
\overline{\Opt}[\Id-H^*A]\leq\gamma
\ee{oldverif}
on the matrix variable $H$, where $\overline{\Opt}$ is given by (\ref{overline}) with $\Z^*_k=\U_k$, is a verifiable sufficient condition for the validity of (\ref{eq500}).
\par
Unfortunately, the sufficient condition in (\ref{oldverif}) is really poor.
Note that in the context of Compressive Sensing, the only interesting case is the one of $\gamma<1$. We are about to show that (\ref{oldverif}) can certify the validity of (\ref{eq500}) with $\gamma<1$ only in  extreme cases:
\begin{proposition}\label{pbadnews}
In the situation in question, it holds
\begin{equation}\label{badnews}
\overline{\Opt}[\Id-H^*A]\geq \min\left[2s\sqrt{{\dim(\Ker (A))\over pq}},\sqrt{\dim(\Ker (A))}\right].
\end{equation}
In other words, when the dimension of $\Ker (A)$ is of order of the dimension $pq$ of $X$, our verifiable sufficient condition for the validity of the nuclear norm recovery can certify $s$-goodness of $A$ only for $s=O(1)$.
\end{proposition}
The proof of the proposition is provided in the appendix.
\par
In fact, the tractable convex relaxations $\Z_k\mapsto \U_k$ underlying the condition (\ref{oldverif}) can be improved, resulting in a weakened form of (\ref{oldverif}). Specifically, let $M'[h,w]=h\otimes w\in \bbr^{p^2\times q^2}$. Obviously, if $M[h,w]=h^T\otimes w$, then $M'[h,w]=\cM'(M[h,w])$ for certain linear map $\cM':\,\bbr^{pq\times pq}\to \bbr^{p^2\times q^2}$ induced by an appropriate $p^2\times q^2$ rearrangement of the entries in a $pq\times pq$ matrix. Using again the fact that the singular values of $h\otimes w$ are the pairwise products of the singular values of $h$ and $w$ we conclude that the set
\[
\V_k=\{\cU%
\in \U_k:\;\|\cM'(\cU)\|\le k,\,\|\cM'(\cU)\|_{2,2}\le 1\}
\]
contains $\Z_k$. Further, let $f(h)=\Col(h^T)=[\Row_1(h)^T;...;\Row_p(h)^T]$, where $\Row_i(h)$ is the $i$-th row of $h$; let also $g(w)=\Col(w)=[\Col_1(w),...,\Col_q(w)]$, where $\Col_j(w)$ is the $j$-th column of $w$. Let us denote $M''[h,w]=f(h)g(w)^T$. It is immediately seen that $M''[h,w]=\cM''(M[h,w])$ for certain linear mapping $\cM'':\,\bbr^{pq\times pq}\to \bbr^{pq\times pq}$ induced by an appropriate permutation of the entries in a $pq\times pq$ matrix. Furthermore,
\[
\max_{h,w}\left\{\|M''[h,w]\|\,\,\,\bigg|\begin{array}{l}\;h\in\H_k=\{h\in\bbr^{p\times q}:\|h\|\le k,\|h\|_{2,2}\le 1\},\\
w\in\W=\{w\in\bbr^{p\times q}:\|w\|\le 1\}\\
\end{array}\right\}
\leq \sqrt{k}.
\]
\begin{quote}{\small
Indeed, since the mapping $[h,w]\to M''[h,w]$ is bilinear, and $\H_k$ and $\W$ are convex, it suffices to verify that whenever $h$ is an extreme point of $\H_k$ and $w$ is the extreme point of $\W$, we have $\|M''[h,w]\|\leq\sqrt{k}$. To this end note that an extreme point $h$ of $\H_k$ is of the form $\sum_{\ell=1}^kh^\ell=\sum_{\ell=1}^ka^\ell [b^\ell]^T$ with unit mutually orthogonal vectors $a^\ell\in\bbr^p$ and unit mutually orthogonal vectors $b^\ell\in\bbr^q$, while an extreme point $w$ of $\W$ is of the form $a b^T$ with unit vectors $a\in \bbr^p$ and $b\in\bbr^q$. By the definition of  $M''[h,w]$ we get
\[
M''[h,w]= f(h)g(w)^T=\left[\sum_{\ell=1}^kf(h^\ell)\right] g^T,
\]
where $f(h^\ell)=[a^\ell_1b^\ell;...;a^\ell_pb^\ell]$ and $g(w)=[ab_1;...;ab_q]$. We see that $f(h^\ell)$ are mutually orthogonal unit $pq$-dimensional vectors, and $g(w)$ is a unit $pq$-dimensional vector, so that $M''$ is a rank 1 matrix with the maximal singular value $\sqrt{k}$. Thus, $\|M''[h,w]\|\leq\sqrt{k}$, as claimed.
}\epr
\end{quote}
Let now
\bse
\Z^*_k&=&\left\{\cU\in  \V_k:\,\|\cM''(\cU)\|\le \sqrt{k}\right\}\\
&=&\left\{\cU\in \bbr^{pq\times pq}\;\left|\begin{array}{l}\|\cU\|\le 1,\,\|\cU\|_{2,2}\le 1;\\
\|\cM'(\cU)\|\le k,\,\|\cM'(\cU)\|_{2,2}\le 1;\\
\,\|\cM''(\cU)\|\le \sqrt{k}\end{array}\right.
\right\}.
\ese
We conclude from the above that $\Z^*_k$ is a convex set which cointains $\Z_k$ and is contained in $\U_k$. Thus the function
\[
\Opt^*[W]=\max_{\cU,\cV}\left\{\sum_{i,j=1}^{pq}(\Theta[W])_{ij}[\cU_{ij}+\cV_{ij}]: \;\cU\in \Z^*_{s},\,\cV\in \Z^*_{2s}\right\}
\]
is an efficiently computable  upper bound on ${\Opt}[W]$ which is $\leq\overline{\Omega}[W]$. We arrive at
\begin{proposition}\label{proplast} If $\widehat{\beta}[H]$ is an efficiently computable convex in $H$ upper bound on $\beta[H]$, then solvability of the system of efficiently computable convex constraints
\[
\Opt^*[\Id-H^*A]\leq\gamma,\;\;\widehat{\beta}[H]\leq \beta
\]
in matrix variable $H$, $\beta$ and $\gamma$ being parameters, is a verifiable sufficient condition for the validity of {\rm (\ref{elln})}.
\end{proposition}
\paragraph{Discussion}
Though  we do not know at this point whether the verifiable condition of Proposition \ref{proplast} for the validity of (\ref{elln}) still obeys the ``limits of performance'' of Proposition\ref{pbadnews}, numerical experiments show that the computationally tractable set $\Z^*_k$ is a proper subset of $\U_k$. Moreover,
no one of the constraints in the description of $\Z^*_k$ is redundant, meaning that dropping any one of them can strictly increase
the solution set.
At the moment $Z^*_k$ is the smallest computationally tractable relaxation of $\Z_k$ known to us.
\par
Note that the verifiable sufficient conditions we have presented for Examples I.a-b can certify $s$-goodness of a sensing map for $s$ as large as $O(\sqrt{m})$ (for details, see \cite{JuN1,Judetal}).
On the other hand, it should be noted that not so good (or even really bad, as in the case of verifiable condition (\ref{oldverif})) ``limits of performance'' of our verifiable sufficient conditions for the validity of sparsity/low rank oriented recovery do not automatically mean that these conditions are of small or no interest. The point is that the recovery in question makes perfect sense also in the case when the observations are not ``deficient'', that is, the sensing map $A$ has the trivial kernel $\{0\}$. Since such an $A$ can be poorly conditioned, the question of how well can we recover (representations of)  sparse/low rank signals remains nontrivial and important even in this ``fully observable'' case. Here our verifiable sufficient conditions, at least in the case of $B=\Id$, can certify the validity of, say, $\bC^+_s(1/2,\beta,\phi)$ for all $s$ and all $\beta\geq\beta(s)$, with some efficiently computable $\beta(s)$, and thus we can efficiently upper-bound the recovery errors as functions of signal's sparsity.

\appendix
\section{Proofs}

\subsection{Proof of Proposition \ref{proprec}}
Let $P\in\cP_s$ be such that $\|(I-P)Bx\|\leq\delta_x$, and let $\bar{w}=PBx$, $\widetilde{w}=Bx-PBx$, $z=\widehat{x}-x$. We have
 \begin{equation}\label{summary}
 \begin{array}{ll}
 (a)&P\bar{w}=P^2Bx=PBx=\bar{w}~~~\hbox{\ [since all $P\in\cP$ are projectors]}\\
 (b)&\|B\widehat{x}\|\leq \Opt+\delta\leq \|Bx\|+\delta ~~~\hbox{\ [since clearly $\Opt\leq\|Bx\|$]}\\
 (c)&\phi(Az)\leq\phi(A\widehat{x}-y)+\phi(y-Ax)\leq \epsilon+\delta_\phi+ \epsilon=2\epsilon+\delta_\phi\\
 \end{array}
 \end{equation}
We now have
\bse
 \|Bx\|+\delta&\geq& \|B\widehat{x}\|=\|Bx+Bz\|=\|\bar{w}+\widetilde{w}+Bz\| ~~~\hbox{\ [by (\ref{summary}.$b$)]} \\
&\geq& \|\bar{w}+Bz\|-\|\widetilde{w}\|~~~\hbox{\ [by the triangle inequality]}\\
 &\geq& -\|\widetilde{w}\|+\|\bar{w}\|+\|\overline{P}Bz\|-\|PBz\|~~~\hbox{\ [by (\ref{summary}.$a$) and Lemma \ref{lem1}]}\\
 &\geq&  -\|\widetilde{w}\|+\|\bar{w}\| +(1-\gamma)\|Bz\|-\beta\phi(Az)~~~\hbox{\ [by $\bC_s(\gamma,\beta;\phi)$]}.
 \ese
 Thus,
\bse (1-\gamma)\|Bz\| &\leq& \|Bx\|+\delta-\|\bar{w}\|+\|\widetilde{w}\|+\beta\phi(Az)
 \leq \delta+2\|\widetilde{w}\|+\beta[2\epsilon+\delta_\phi]\\
&&~~~\hbox{\ [by (\ref{summary}.$c$) and the triangle inequality; note that $Bx=\bar{w}+\widetilde{w}$]}\\
 &\leq& \delta+2\delta_x+\beta[2\epsilon+\delta_\phi]~~~\hbox{\ [since $\|\widetilde{w}\|=\|Bx-PBx\|\leq\delta_x$]},
\ese
and (\ref{errorbound}) follows. \epr
\subsection{Proof of Proposition \ref{penrec}}
Let $x$, $P\in\cP_s$, $\widehat{x}$ satisfy the premise of the proposition, and let $z=\widehat{x}-x$, $\xi=y-Ax$. We have
\bse
\|B\widehat{x}\|+\lambda\phi(A\widehat{x}-y)&\leq& \|Bx\|+\lambda\phi(Ax-y)+\delta,\\
\phi(A\widehat{x}-y)&=&\phi(Az+Ax-y)\geq \phi(Az)-\phi(\xi),
\ese
and we conclude that
\be
\|B\widehat{x}\|+\lambda\phi(Az)\leq \|Bx\|+2\lambda\phi(\xi)+\delta.
\ee{*}
Further,
\bse
\|B\widehat{x}\|&=&\|Bx + Bz\|=\|PBx+(I-P)Bx+Bz\|\\
&\geq& \|PBx+Bz\|-\|(I-P)Bx\|\;\;~\hbox{\ [by the triangle inequality]}\\
&\geq& -\|(I-P)Bx\|+\|PBx\|+\|\overline{P}Bz\|-\|PBz\|
\ese
(the latter inequality is (\ref{eq1}) as applied with $w=PBx$; note that $Pw=w$ due to $P^2=P$).
Thus, when applying condition $\bC_s(\gamma,\beta;\phi)$, we obtain
\bse
\|B\widehat{x}\|
&\geq& -\delta_x+\|PBx\|+(1-\gamma)\|Bz\|-\beta\phi(Az),
\ese
and therefore
\bse
-\delta_x+\|PBx\|+(1-\gamma)\|Bz\|-\beta\phi(Az)+\lambda\phi(Az) &\leq& \|B\widehat{x}\|+\lambda\phi(Az)\\
&\leq& \|Bx\|+2\lambda\phi(\xi)+\delta
\ese
where the concluding $\leq$ is due to (\ref{*}). We come to
\bse
(1-\gamma)\|Bz\|&\leq& \|Bx\|-\|PBx\|+\delta_x+\underbrace{[\beta-\lambda]}_{\leq0}\phi(Az)+2\lambda\phi(\xi)+\delta\\
&\leq&2\delta_x+2\lambda\phi(\xi)+\delta \hbox{\ [the triangle inequality]}
\ese
The resulting relation is nothing but (\ref{imperfect}). \epr
\subsection{Proof of Proposition \ref{pbadnews}}
Let us treat $pq\times pq$ matrices $\cA$ as $q\times p$ block matrices  with $p\times q$ blocks $\cA^{\mu\nu}$, $1\leq\mu\leq q,1\leq \nu\leq p$, so that
for $\cA=w^T\otimes z$ we have $\cA^{\mu\nu}=w_{\nu\mu}z$. Now, let $X=\bbr^{p\times q}$, and let $L(X,X)$  be the family of all linear maps $z\to Rz:X\to X$. A map $R\in L(X,X)$  can be identified with tensor $R_{ijk\ell}$, $1\leq i,k\leq p$, $1\leq j,\ell\leq q$, according to $[Rz]_{ij}=\sum_{k,\ell}R_{ijk\ell}z_{k\ell}$. By definition of $\Theta[\cdot]$, for every $R\in L(X,X)$ and all $z,w\in X$ we should have
\[
\Tr(\Theta[R][w^T\otimes z])=\Tr([Rz]w^T)=\sum_{i,j,k,\ell} R_{ijk\ell}w_{ij}z_{k\ell},
\]
so that $\Theta[R]$ is $pq\times pq$ matrix with the blocks $(\Theta[R])^{\mu\nu}=R^{\mu\nu}:=[R_{\nu\mu k\ell}]_{{1\leq k\leq p,\atop 1\leq \ell\le q}}\in\bbr^{p\times q}$. In particular, if $\cS$ is a $pq\times pq$ matrix such that all blocks $\cS^{\mu\nu}$ belong to the kernel of $R$, then $\Tr(\Theta[R]\cS)=\sum_{\mu,\nu}\Tr(R^{\mu\nu}[\cS^{\mu\nu}]^T)=0$. \par
Now let $\cE$ be the $pq\times pq$ matrix with the blocks $\cE^{\mu\nu}$ possessing each exactly one nonzero entry, equal to 1, in the position $\nu\mu$, so that $\cE^{\mu\nu}$ are just the standard basic orths in $\bbr^{p\times q}$. Denoting by $\cW^{\mu\nu}$, $\cH^{\mu\nu}$ the blocks in the $pq\times pq$ matrices $\Theta[W]$ and $\Theta[H^*A]$, respectively,  we have for all $z,w\in X$:
\[
\Tr(\cE[w^T\otimes z]^T)=\Tr(zw^T) =\Tr([(W+H^*A)z]w^T)=\Tr(\Theta[W][w^T\otimes z]^T)+\Tr(\Theta[H^*A][w^T\otimes z]^T),
\]
whence $\cE^{\mu\nu}=\cW^{\mu\nu}+\cH^{\mu\nu}$ for all $\mu,\nu$.

Consider any $\cZ\in\Xi$ where $\Xi$ is the set of all $pq\times pq$ matrices $\cZ$ of nuclear norm $\leq 1$ and such that all blocks $\cZ^{\mu\nu}$ belong to the kernel of the sensor map $A$. Then $\cZ^{\mu\nu}$ belong to the kernel of $H^*A$, so that
by the above $\Tr(\cH^{\mu\nu}[\cZ^{\mu\nu}]^T)=0$ for all $\mu,\nu$, whence $\Tr(\cE\cZ^T)=\Tr(\Theta[W]\cZ^T)$.
\par
We clearly have
\bse
\overline{\Opt}[W]\geq \underline{\Opt}[W]&:=&\max_\cZ\left\{\Tr(\Theta[W]\cZ^T):\cZ\in\bbr^{pq\times pq},\|\cZ\|\leq 2s,\|\cZ\|_{2,2}\leq1,\,\cZ^{\mu\nu}\in\Ker (A)\,\forall \mu,\nu
\right\}\\
&=&\max_{\cZ}\left\{\Tr(\cE\cZ^T): \cZ\in\bbr^{pq\times pq},\|\cZ\|\leq 2s,\|\cZ\|_{2,2}\leq1,\cZ^{\mu\nu}\in\Ker (A)\,\forall \mu,\nu\right\}.
\ese
Passing to the dual problem, we get
\be
\underline{\Opt}[W]=\min\limits_{\cU,\cV}\left\{2s\|\cU\|_{2,2}+\|\cV\|:
\cU,\cV\in \bbr^{pq\times pq}, \cU^{\mu\nu}+\cV^{\mu\nu}-\cE^{\mu\nu}\in(\Ker (A))^\perp\,\forall \mu,\nu\right\}.
\ee{eqno}
Now let $\cU,\cV$ be a feasible solution to the latter problem. Denoting by $\|\cdot\|_F$ the Frobenius norm on $\bbr^{pq\times pq}$, we have
\begin{equation}\label{eq510}
\|\cU\|_F^2\leq pq\|\cU\|_{2,2}^2,\,\,\,\|\cV\|_F^2\leq \|\cV\|^2.
\end{equation}
On the other hand, $(\cU+\cV)^{\mu\nu}-\cE^{\mu\nu}\in (\Ker (A))^\perp\,\,\forall\mu,\nu
$
implies that
\be
\|\cU+\cV\|_F^2\geq \sum_{\mu,\nu} \dist^2(\cE^{\mu\nu},(\Ker (A))^\perp),
\ee{**}
where the distance is measured in the Frobenius norm. Let us treat $\bbr^{p\times q}$ as a $pq$-dimensional Euclidean space $\bbr^{M}$, $M=pq$,  equipped with the standard Euclidean structure, and let $d$ be the dimension of $\Ker (A)$. Observe that the right hand side in $(\ref{**})$ is sum, over all basic orths $\cE^{\mu\nu}$ of $\bbr^M$,  of squared distances from the basic orths to a linear subspace $L$ of $\bbr^M$ of codimension $d$, namely, $L=(\Ker (A))^\perp$.  Note that such a sum is at least $d$. Indeed,
otherwise we could approximate the unit $M\times M$ matrix by with a matrix of rank $\leq M-d$ (specifically, a matrix with columns from $L$) within an accuracy, measured  in the Frobenius norm, which is better than $\sqrt{d}$, which is impossible. Thus, the right hand side in $(\ref{**})$ is at least $d=\dim(\Ker (A))$, and we arrive at the inequality
\[
\|\cU+\cV\|_F^2\geq d:=\dim(\Ker (A))~~\Rightarrow~~ \|\cU\|_F+\|\cV\|_F\geq\sqrt{d}.
\]
Invoking (\ref{eq510}) and recalling that $\cU,\cV$ is an arbitrary feasible solution to (\ref{eqno}), we arrive at
\bse
\overline{\Opt}[W]&\geq&\underline{\Opt}[W]\geq \min_{a,b\geq0}\left\{{2s\over \sqrt{pq}}a+b:a+b\geq\sqrt{\dim\Ker (A)}\right\}\\
&=&\min\left[2s\sqrt{{\dim(\Ker (A)) \over pq}},\sqrt{\dim(\Ker (A))}\right].
\ese
\epr


\begin{thebibliography}{c}

\bibitem{CandesRecht2009}
Candes, E., Recht, B., ``Exact matrix completion via convex optimization" -- {\sl Foundations of Computational Mathematics} {\bf 9:6} (2009), 717--772.
\bibitem{CandesRecht2012}
Candes, E., Recht, B., ``Simple bounds for low-complexity model reconstruction" -- {\sl Mathematical Programming Ser. A} {\bf 0025-5610} (2012), 1--13.
\bibitem{CandesTao2009}
Candes, E., Tao, T., ``The power of convex relaxation: Near-optimal matrix completion" --{ \sl IEEE Trans. Inform. Theory} {\bf 56:5} (2009), 2053--2080.
\bibitem{NP1}
Donoho, D., Huo, X., ``Uncertainty principles and ideal atomic decomposition'' --- {\sl IEEE Trans. Inf. Theory}
{\bf 47:7} (2001), 2845--2862.
\bibitem{JuN1}
Juditsky, A., Nemirovski, A. ``On Verifiable Sufficient Conditions for Sparse Signal Recovery via $L_1$ Minimization''  -- {\sl Mathematical Programming Ser. B} {\bf 127} (2011), 57--88.
\bibitem{JuKKN2}
Juditsky, A., K{\i}l{\i}n\c{c}  Karzan, F., Nemirovski, A., ``Verifiable conditions of $\ell_1$-recovery of sparse signals with sign restrictions'' --  {\sl Mathematical Programming Ser. B} {\bf 127} (2011), 89--122.
\bibitem{JuN3}
Juditsky, A., Nemirovski, A. ``Accuracy guarantees for $\ell_1$ recovery'' –--  {\sl IEEE Transactions on Information Theory} {\bf 57:12} (2011), 7818--7839.
\bibitem{Judetal}
Juditsky, A., K{\i}l{\i}n\c{c}  Karzan, F., Nemirovski, A., Polyak, B. (2011), ``Accuracy Guarantees for $\ell_1$ recovery of block-sparse signals'' -- submitted to {\sl Annals of Statistics},\\
E-print: http://www.arxiv.org/PS\_cache/arxiv/pdf/1111/1111.2546v1.pdf\\
For abbreviated version, see\\
Juditsky, A., K{\i}l{\i}n\c{c} Karzan, F., Nemirovski, A., Polyak, B.T., ``On the accuracy of $\ell_1$ filtering of signals with block-sparse structure'' -- in: J. Shawe-Taylor, R.S. Zemel, P. Bartlett, F. Pereira, K.Q. Weinberger, Eds. {\sl Advances in Neural Information Processing Systems 24} (2011), 1666-1674.
\bibitem{Tuncel}
Kong, L., Tuncel, L., Xiu, N.,  ``$s$-goodness for low-rank matrix recovery, translated from sparse signal recovery''\\
E-print: www.optimization-online.org/DB\_FILE/2011/06/3063.pdf
\bibitem{OymakFazel2011}
Oymak, S., Mohan, K., Fazel, M., Hassibi, B., ``A Simplified Approach to Recovery Conditions for Low-rank Matrices'' -- {\sl Proc. Intl. Sympo. Information Theory (ISIT)}, Aug 2011.
\bibitem{}
Raghunandan H. Keshavan, Sewoong Oh, and Andrea Montanari. Matrix completion from a few entries.
In ISITÕ09: Proceedings of the 2009 IEEE international conference on Symposium on Information Theory,
pages 324Ð328, Piscataway, NJ, USA, 2009. IEEE Press. ISBN 978-1-4244-4312-3.
\bibitem{RechtFazel}
Recht, B., Fazel, M., Parrilo, P., ``Guaranteed Minimum-Rank Solutions of Linear Matrix Equations via Nuclear Norm Minimization" -- {\sl SIAM Review} {\bf 52:3} (2010), 471--501.
\bibitem{NP2}
 Zhang, Y., {\sl A simple proof for recoverability of $\ell_1$-minimization: go over or under?} --  Rice CAAM
Department Technical Report TR05-09 (2005).
\end{thebibliography}
\end{document}